\documentclass{amsart}
\usepackage{color}
\usepackage[margin=1in]{geometry}
\usepackage{amssymb}
\usepackage{graphicx}
 \usepackage{epstopdf}
 \usepackage{tikz-cd}
\usepackage{hyperref}
\newtheorem{theorem}{Theorem}

\theoremstyle{definition}

\theoremstyle{remark}

\numberwithin{equation}{section}
\begin{document}
\title{A Note on Iterated Maps of the  Unit Sphere}


\author{CHAITANYA GOPALAKRISHNA}
\address{Indian Statistical Institute, Stat-Math. Unit, R V College Post, Bengaluru 560059, India}
\email{cberbalaje@gmail.com, chaitanya\_vs@isibang.ac.in}

\thanks{The author is
    supported
    by the National Board for Higher Mathematics, India through No:
    0204/3/2021/R\&D-II/7389. The author is very grateful to 
    his mentor Professor B. V. Rajarama Bhat for useful discussions.}





\begin{abstract}
Let $\mathcal{C}(S^{m})$ denote the set of continuous maps from the unit sphere $S^{m}$ in $\mathbb{R}^{m+1}$ into itself endowed with the supremum norm. 
We prove that 
the set $\{f^n: f\in \mathcal{C}(S^{m})~\text{and}~n\ge 2\}$ of iterated maps is not dense in $\mathcal{C}(S^{m})$. This, in particular, proves that the periodic points of the iteration operator of order $n$ are not dense in $\mathcal{C}(S^m)$ for all $n\ge 2$, providing an alternative proof of the result that these operators are not Devaney chaotic on $\mathcal{C}(S^m)$ proved in [M. Veerapazham,  C. Gopalakrishna, W. Zhang, 
Dynamics of the iteration operator on the space of continuous self-maps,
{\it Proc. Amer. Math. Soc.},  149(1) (2021), 217--229].
\end{abstract}

\maketitle


\noindent
{\it Iteration}, which refers to repeating the same action, is not only a crucial operation in 
contemporary industrial production but also a typical loop program in computer algorithms.  
Many authors have explored numerous interesting and complex characteristics of this operation through discussions on various aspects such as dynamical systems (\cite{Li,Veer2020}), iterative roots (\cite{BG,Blokh,rice1980,wzhang1997}), and solutions to iterative equations (\cite{Baron-Jarczyk,kuczma1990,Targonski1981,zdun-soalrz}).
An {\it iterated map} on a non-empty set $X$ is simply  a self-map on it of the form $f^n$ for some self-map $f$ on $X$ and an integer $n\ge2$, with $f^k$ denoting the $k$-th order iterate of $f$ defined recursively by $f^0={\rm id}$, the identity map on $X$, and $f^k=f\circ f^{k-1}$ for $k\ge 1$.  
Let $\mathcal{C}(X)$ denote the set of continuous maps from  a locally compact Hausdorff space $X$ into itself in the compact-open topology  and $\mathcal{W}(X):=\cup_{n=2}^{\infty}\mathcal{W}(n;X)$, the set of iterated maps in $\mathcal{C}(X)$, where $\mathcal{W}(n;X):=\{f^n:f\in \mathcal{C}(X)\}$ for all $n\ge 2$. Then the fact that even complex quadratic polynomials are not iterated maps on $\mathbb{C}$, as shown in \cite[Theorem 1]{rice1980}, prompts us to ask the natural and interesting question: {\it how large is $\mathcal{W}(X)$  in $\mathcal{C}(X)$?}
Many researchers have investigated this problem 
from topological and measure theoretic perspectives, and 
the currently known findings in this regard are summarized in the following.

\begin{theorem}\label{T1-nowheredense}
	\begin{enumerate}
		\item[\bf (i)] 		 {\rm (Humke and Laczkovich \cite{Humke-Laczkovich1989,Humke-Laczkovich})}  
		The complement of $\mathcal{W}(2;[0,1])$ is
		dense in $\mathcal{C}([0,1])$;  $\mathcal{W}(2;[0,1])$ is not dense in $\mathcal{C}([0,1])$; and	$\mathcal{W}(n;[0,1])$ is an analytic non-Borel subset of $\mathcal{C}([0,1])$ for $n\ge 2$.
		
		
		
		\item[\bf (ii)] 	
		{\rm  (Simon \cite{simon1989,simon1990,simon1991a,simon1991b})} $\mathcal{W}([0,1])$ is of first category and of zero Wiener measure;  $\mathcal{W}(2;[0,1])$ is nowhere dense in $\mathcal{C}([0,1])$; and
		$\mathcal{W}([0,1])$ is not  dense in $\mathcal{C}([0,1])$. 
		
		\item[\bf (iii)]  	 {\rm  (Blokh \cite{Blokh})} $\mathcal{W}([0,1])$ is nowhere dense  in $\mathcal{C}([0,1])$.
		
		\item[\bf (iv)]  	  {\rm (Bhat and Gopalakrishna \cite{BG,BG-single-preprint})} 
		Every non-empty open set of $\mathcal{C}(X)$ contains a map that is not an iterate of
		even a discontinuous map on $X$ in the following cases: {\bf (a)}  $X=[0,1]^{m+1}$ for $m\ge 0$; {\bf (b)} $X=\mathbb{R}^{m+1}$ for $m\ge 0$; and {\bf (c)} $X=S^1$.  In particular, 	the complement of $\mathcal{W}(X)$ is
		dense in $\mathcal{C}(X)$ in these cases.
	\end{enumerate}
\end{theorem}

Our aim here is to prove the following result, which does not seem to have been proven earlier by anyone. 
\begin{theorem} \label{Not-dense}
	$\mathcal{W}(S^{m})$ is not dense in $\mathcal{C}(S^{m})$ for all $m\ge 1$. 
\end{theorem}

The key to the proof is to use the well-known concept of degree of maps in $\mathcal{C}(S^{m})$, which was first introduced by Brouwer \cite{Brouwer1911}. Each map $f\in \mathcal{C}(S^m)$, as shown in \cite[p.134]{Hatcher}, induces a homomorphism $f_*$ from the singular homology group $H_m(S^m)=\mathbb{Z}$ of $S^m$ into itself, which has the form $f_*(a)=ka$ for some $k\in \mathbb{Z}$. The integer $k$ is called the {\it degree} of $f$, with the notation ${\rm deg}(f)$.
Intuitively, ${\rm deg}(f)$ counts how many times $f$ wraps $S^m$ around itself, with the sign indicating whether $f$ preserves orientation or not. 
Interestingly, unlike the proofs of the results in Theorem \ref{T1-nowheredense}, which are both technical and involved, the proof of Theorem \ref{Not-dense} is  simple and does not require any complex machinery, as shown below.  

%
%

%

\begin{proof}[Proof of Theorem 2]
	Since the degree map on $\mathcal{C}(S^m)$ satisfies ${\rm deg}(f\circ g)={\rm deg}(f)\cdot {\rm deg}(g)$ as shown in \cite[p. 134]{Hatcher}, we see that the degree of each $f\in \mathcal{W}(S^m)$ has the form $k^n$ for some $k\in \mathbb{Z}$ and $n\ge 2$. 
	Consider an $f_0\in \mathcal{C}(S^{m})$ such that ${\rm deg}(f_0)=2$, which exists due to Example 2.31 in \cite{Hatcher}, and let $B_1(f_0)$ be the unit ball in $\mathcal{C}(S^m)$ centered at $f_0$. Then, for each $g\in B_1(f_0)$, the map $H_g:S^m\times [0,1]\to S^m$ defined by
	\begin{eqnarray*}
		H_g(x,t)= \frac{(1-t)f_0(x)+tg(x)}{\|(1-t)f_0(x)+tg(x)\|_2},\quad \forall (x,t)\in S^m\times [0,1]
	\end{eqnarray*}
	is a homotopy between $f_0$ and $g$, implying that ${\rm \deg}(g)={\rm deg}(f_0)$  by Property (d) in \cite[p. 134]{Hatcher}, which states that homotopic maps have the same degree. 
	Therefore, as $2\ne k^n$ for all $k\in \mathbb{Z}$ and $n\ge 2$, we have $B_1(f_0)\cap \mathcal{W}(S^{m})=\emptyset$. 
\end{proof}

The choice of $f_0$ with ${\rm deg}(f_0)=2$ in the proof of Theorem \ref{Not-dense} is not special. Indeed, the argument works even for an $f_0$ with ${\rm deg}(f_0)=l$ for any integer $l$ that is not a perfect power, demonstrating that there are infinitely many open balls in $\mathcal{C}(S^{m})$ each element of which is not an iterated map  in $\mathcal{C}(S^{m})$. 
Furthermore, although the result in Theorem \ref{Not-dense} appears to have been obtained for free by using only well-known results in singular homology,
we believe it is still worth recording because it is relevant to the 
{\it iterative root problem} (i.e., the problem of finding a self-map $f$ satisfying $f^n=F$ for a given self-map $F$ and an integer $n\ge 2$) 
 rooted in Babbage's classic work \cite{Babbage1815} and is complimented by similar results in Theorem \ref{T1-nowheredense}.
Additionally, whereas all known results particularly for spheres (see \cite{Jarczyk2003,solarz2002,solarz2003,solarz2006,solarz2012,Zdun2000,zdun2008} and the references therein) are either limited to $S^1$ or are based on the concept of {\it rotation number}, our result deals with $S^m$ for all $m\ge1$ and employs a different tool of {\it degree} of maps.

 On the other hand, Theorem \ref{Not-dense} is significant in the theory of dynamical systems as well for the following reason: As in \cite{Devaney},
 a discrete
 semi-dynamical system $(X,f)$, where $X$ is a metric space equipped with the metric $d$,
 is said to be {\it topologically transitive} if for every pair of open sets $U, V$ in $X$ there exist $x \in U$ and $n \in \mathbb{N}$
 such that $f^n(x)\in V$.
 $f$ is said to be {\it sensitively dependent on initial conditions} if there exists $\delta >0$ such that
 for every $x \in X$ and every $\epsilon>0$
 there exist $y \in X$ and $n \in \mathbb{N}$ such that $d(x,y)<\epsilon$ and $d\big(f^n(x), f^n(y)\big)>\delta$.
 We say that $f$ is {\it chaotic} in Devaney's sense
 if {\bf (i)} the set of periodic points of $f$ is dense in $X$, {\bf (ii)} $f$ is topologically transitive, and
 {\bf (iii)} $f$ exhibits sensitive dependence on initial conditions. The authors of \cite{Veer2020} investigated the dynamics of {\it iteration operators} $\mathcal{J}_n:=f^n$ of orders $n\ge 1$ and proved that these operators are neither topologically transitive nor exhibit sensitive dependence on initial conditions, establishing that they are not Devaney chaotic on $\mathcal{C}(X)$ for compact metric spaces $X$ (see \cite[Theorem 5.1]{Veer2020}).
 However, since no maps in $B_1(f_0)$ are iterated maps, as shown in the proof, Theorem \ref{Not-dense} provides an alternative proof of this result for the spheres $S^m$ and $n\ge 2$ by proving that the periodic points of $\mathcal{J}_n$ are not dense in $\mathcal{C}(S^m)$.

Besides, although the result in Theorem \ref{Not-dense} is given for the spheres $S^m$, we can anticipate a similar result and its implications (for the iterative root problem and the dynamics of iteration operators) for  connected closed orientable $k$-manifolds, in which context the degree for continuous maps is defined (see Exercise 7 in \cite[p.258]{Hatcher}). Moreover, we believe that Theorem \ref{Not-dense} can be presented in the first course on homology theory immediately after introducing singular homology groups,  as one of its typical applications to the theory of iterative roots, providing graduate students with an early opportunity to explore iteration theory and consider it as a potential area 
of
research.
Finally, we conclude this note with the following open problem for future discussion, which is inspired by result {\bf (iii)} of Theorem \ref{T1-nowheredense}: {\it Is $\mathcal{W}(S^m)$ nowhere dense in $\mathcal{C}(S^m)$ for $m\ge 1$?}

\bibliographystyle{amsplain}

\end{document}